
\input gtmacros

\input gtmonout
\volumenumber{2}
\volumeyear{1999}
\volumename{Proceedings of the Kirbyfest}
\pagenumbers{563}{569}
\papernumber{26}
\received{28 July 1998}
\published{21 November 1999}

\def\Z{{\Bbb Z}}

\def\R{{\Bbb R}}
\def\cy#1{\Z/{#1}\Z}

\reflist

\key{bBo}
{\bf A~Borel}, {\it Seminar on {T}ransformation {G}roups}, Annals of
Math. Studies, vol.~46, Princeton Univ. Press, Princeton (1960)

\key{bBr}
{\bf G\,E Bredon}, {\it Introduction to {C}ompact {T}ransformation
{G}roups}, Academic\break Press, New York (1972)

\key{bDW}
{\bf D~van Danzig}, {\bf B\,L~van~der Waerden}, {\it \"{U}ber
metrische homogene R\"aume}, Abh.  Math. Sem. Hamburg, {6} (1928)
374--376

\key{bFT}
{\bf M\,H Freedman}, {\bf L\,R Taylor}, {\it A universal smoothing of 
four{--}space}, J. Differential Geometry, {24} (1986) 69--78

\key{bKN}
{\bf S~Kobayashi}, {\bf K~Nomizu}, {\it Foundations of
{D}ifferential {G}eometry, Volume I}, John Wiley, New York (1963)

\key{bMS}
{\bf S\,B Myers}, {\bf N\,E Steenrod}, {\sl The group of isometries of
a {R}iemannian manifold}, Ann. of Math. {40} (1939)
400--416

\key{bOl}
{\bf R Oliver}, {\it A proof of the Conner conjecture}, Ann.of Math. {103}
(1976) 637--644

\key{bCT}
{\bf C\,H Taubes}, {\it Gauge theory on asymptotically periodic
4--manifolds}, J.  Differential Geom. {25} (1987) 363--430

\key{bT}
{\bf L\,R Taylor}, {\it An invariant of smooth $4${--}manifolds},
Geometry and Topology, {1} (1997) 71--89

\endreflist

\def \fbic{Theorem\ 1.1} 
\def \fsoa{Theorem\ 2.1} 
\def \fbst{Theorem\ 3.1} 
\def \fiso{Lemma\ 3.2} 
\def \fdso{Lemma\ 3.3} 

\newbox\nbox
\def\Bar#1{\setbox\nbox = \hbox{$#1$}
\kern .08\wd\nbox{\overline{\hbox to .84\wd\nbox{\vphantom {$#1$}\hss}}}
\kern .08\wd\nbox\kern -\wd\nbox\box\nbox}

\let\disjointunion\coprod

\def\newT#1{\relax}
\def\newL#1{\relax}

\let\newsec\section

\title{Smooth Euclidean 4{--}spaces with few symmetries}
\asciititle{Smooth Euclidean 4-spaces with few symmetries}
\author{Laurence R Taylor}
\address{Department of Mathematics, University of Notre Dame\\ 
Notre Dame, IN 46556, USA}
\email{taylor.2@nd.edu}

\abstract

We say that a topologically embedded 3--sphere in a smoothing of
Euclidean 4--space is a barrier provided, roughly, no diffeomorphism
of the 4--manifold moves the 3--sphere off itself.  In this paper we
construct infinitely many one parameter families of distinct
smoothings of 4--space with barrier 3--spheres.

The existence of barriers implies, amongst other things, that the
isometry group of these manifolds, in any smooth metric, is finite.
In particular, $S^1$ can not act smoothly and effectively on any
smoothing of 4--space with barrier 3--spheres.

\endabstract
\asciiabstract{We say that a topologically embedded 3-sphere in a 
smoothing of Euclidean 4-space is a barrier provided, roughly, no
diffeomorphism of the 4-manifold moves the 3-sphere off itself.  In
this paper we construct infinitely many one parameter families of
distinct smoothings of 4-space with barrier 3-spheres.

The existence of barriers implies, amongst other things, that the
isometry group of these manifolds, in any smooth metric, is finite.
In particular, S^1 can not act smoothly and effectively on any
smoothing of 4-space with barrier 3-spheres.}


\maketitle

We say that a smooth manifold has {\sl few symmetries\/} provided
that, for every choice of smooth ($C^1$ or better) metric, 
the isometry group for that metric is finite.
Let $E^4$ be a smooth manifold homeomorphic to $\R^4$.
We say that a flat embedding $S^3\subset E$ is a {\sl barrier\/} $S^3$
provided that, given any open set $U\subset E^4$ containing $S^3$
and any smooth embedding $e\colon U\to E^4$, then 
$e(S^3)\cap S^3\neq\emptyset$.
Given a barrier $S^3\subset E^4$, the {\sl inside\/} is the 
component of $E^4-S^3$ whose closure is compact: note that it
is a smoothing of $\R^4$.

\proclaim{Theorem}
Let $E^4$ be a smoothing of $\R^4$ with a barrier $S^3$ 
whose inside does not smoothly 
embed in any integral homology $4${--}sphere.
Then $E$ has few symmetries.\endproc

\noindent{\bf Remarks}\qua
In [\bT], we have constructed many examples of 
smooth $\R^4$'s as in the theorem.
See the discussion leading up to \fbst\ below and the theorem itself.
There definitely are examples for which the isometry group is
not trivial: eg, the end-connected sum of $E$ with itself supports
an involution and any metric can be averaged so as to make the
involution an isometry.

\prf
By \fbic\ below,
the barrier forces the isometry group to be a compact Lie group.
Myers and Steenrod [\bMS] have already proved that the isometry
group is a Lie group, so our contribution is that it must be compact.
By \fsoa\ below, if $S^1$ acts effectively on $E^4$, then
any compact subset of $E^4$ embeds smoothly in an integral
homology $4${--}sphere.
Hence, under our hypotheses, the isometry group is 
a compact Lie group with no $S^1$
subgroups, which implies that it is finite.\endprf

\newsec{Barriers and isometry groups.}
The goal of this section is to prove:
\newT{bic}
\proclaim\fbic
Fix a $C^1$ metric on a smooth $\R^4$ with a barrier $S^3$.
Then the isometry group in this metric is compact.\endproc

\prf Myers and Steenrod [\bMS] have proved that the isometry
group is a Lie group.
Hence it will suffice to show that any infinite set of isometries
has a convergent subsequence.
Let ${\cal I}$ denote any infinite set of isometries.

Fix $\epsilon>0$ so that $C_\epsilon=
\{ x\in E\ \vert\ d(x,S^3)\leq\epsilon\}$
is compact.
If the metric is complete, any finite number will do, but even if
the metric is not complete, there are sufficiently small $\epsilon$
with this property.
Cover $S^3$ with finitely many balls, 
$B(x_i,\epsilon/2)$ with $x_i\in S^3$, $i=1$, \dots, $r$.
Let ${\cal I}_i=\{ f\in{\cal I}\ \vert\ f\bigl(B(x_i,\epsilon/2)\bigr)
\subset E^4-S^3\}$.
If $f$ is an isometry and if $f\notin{\cal I}_i$, equivalently
$f\bigl(B(x_i,\epsilon/2)\bigr)\cap S^3\neq\emptyset$,
then $f\bigl(B(x_i,\epsilon/2)\bigr)\subset C_\epsilon$.

There exists an $i$ such that ${\cal I} - {\cal I}_i$ 
is infinite.
To see this, observe that
$\cap_{i=1}^r{\cal I}_i=\emptyset$ because $S^3$ is a barrier
and $\cup_{i=1}^r\bigl({\cal I}-{\cal I}_i\bigr)=
{\cal I}-\cap_{i=1}^r{\cal I}_i={\cal I}$.
Since ${\cal I}$ is infinite, so is at least one
${\cal I} - {\cal I}_i$.
Let ${\cal I}^0={\cal I} - {\cal I}_i$ for any $i$ such that
${\cal I} - {\cal I}_i$ is infinite.

Pick $5$ points, $y_0$, \dots, $y_4\in {\cal I}^0$ in sufficiently general
position as discussed by Myers and Steenrod in Theorem 3, [\bMS].
Then ${\cal I}^0(y_t)=\{ f(y_t) \ \vert\ \forall f\in{\cal I}^0\ \}
\subset C_\epsilon$ so there is a subset of ${\cal I}^0$,
$f_0$, \dots, $f_r$, \dots so that for each $t$, 
the sequence $f_r(y_t)$ is a Cauchy sequence.
As discussed by Myers and Steenrod [\bMS] at the top of page 406,
it follows from a theorem of van Danzig and van der Waerden 
[\bDW] that there exists an isometry $f$ such that 
a subsequence of the $f_r$ converge to $f$.
(A proof can also be found in [\bKN], especially the proof of Theorem 4.7
starting on page 46.)
Since every infinite subset of the isometry group has a
convergent subsequence, it follows that the isometry group
is compact.\endprf

\newsec{Effective $S^1$ actions on smoothings of $\R^4$'s.}
In this section we will show:
\newT{soa}
\proclaim\fsoa
If $S^1$ acts effectively on a smoothing $E$ of $\R^4$
then any compact smooth submanifold of $E$ embeds
smoothly in an integral homology $4${--}sphere.\endproc

The proof occupies the remainder of this section.
We begin with some general results on $S^1$ actions.
First recall that if $S^1$ acts effectively on
a connected manifold $M$, then the dimension of each
component of $M^{S^1}$ is congruent mod $2$ to the dimension
of $M$.
If $H\subset S^1$ is a proper subgroup, then it is finite cyclic.
If a component of $M^H$ contains a component of $M^{S^1}$
then the dimension of that component of $M^H$ is congruent 
mod $2$ to the dimension of $M$.
To see this, let $C$ denote the component of $M^H$ and 
let $\widehat H$ denote the subgroup for which
$S^1/\widehat H$ acts effectively on $C$.
Then $\bigl(C\bigr)^{S^1/\widehat H}=C\cap M^{S^1}$ so the
dimension of $C$ is congruent mod $2$ to the dimension of a
component of $M^{S^1}$.

Now suppose $M$ is acyclic over the integers.
Then so is $M^{S^1}$: in particular it is non-empty and connected.
For any prime $p$, $M^{\cy p}$ is mod $p$ acyclic, hence also
non-empty and connected.
If the codimension of $M^{S^1}$ in $M$ is $2$, then the
action must be semi-free since $M^{\cy p}=M^{S^1}$ for
all primes $p$ and hence $M^H=M^{S^1}$ for 
any non-trivial, proper subgroup $H\subset S^1$.

The only other case of relevance here is the case where
$M$ is still integrally acyclic, the dimension of $M$
is $4$, and the dimension of $M^{S^1}$ is $0$.
Some of the $M^{\cy p}$ may have dimension $2$, but
there are only finitely many, say $p_1$, \dots, $p_r$.
Let $K(p_i)$ denote the subgroup so that $S^1/K(p_i)$ acts
effectively on $M^{\cy{p_i}}$.
Now $M^{\cy{p_i}}$ is a mod $p_i$ acyclic, non-compact
$2${--}manifold, hence $\R^2$.
Conveniently, $M^{\cy{p_i}}$ remains integrally acyclic, so the $S^1/K(p_i)$ 
action must be semi-free.
Hence either
$M^{K(p_i)}\cap M^{K(p_j)}=M^{S^1}$ or $K(p_i)=K(p_j)$.
By a theorem of Bochner's, the action in a neighborhood of $M^{S^1}$
is linear, so the intersections are transverse as well.

If $M^4$ is actually contractible, then so is the orbit space
$M^\ast$ [\bOl; page 644 Theorem 5].
By [\bBr; page 189, 4.6], $M^\ast$ is a $3${--}manifold, with boundary
if $M^{S^1}=\R^2$ and without boundary if $M^{S^1}$
is a point.
The image of $M^{S^1}$ in $M^\ast$ is the boundary if $M^{S^1}=\R^2$.
Let $P\subset M$ be the set of principal orbits:
$P$ is an open dense set.
If $P^\ast=\pi(P)$, then $P^\ast$ is an open dense set and the map
$\pi\colon P\to P^\ast$ is a submersion. 
Since $M^\ast$ is orientable, for any any embedded 
$S^1\subset P^\ast$, $\pi^{-1}(S^1)$ is a torus, not a Klein bottle:
ie, for any $S^1\subset P^\ast$, the circle action on
$\pi^{-1}(S^1)$ is trivial.

Let us now restrict attention to the case $M=E$ is homeomorphic
to $\R^4$.
Let $U$ be the interior of a topological ball in $E$ which contains
our compact submanifold.
It will suffice to embed $U$ smoothly 
in an integral homology $4${--}sphere.
Let $U^\ast$ denote the image of $U$ in $E^\ast$.

Pick a point $p\in E^{S^1}$ and let $V\in E$ be a 
smooth linear $4${--}ball centered at $p$.
Let $V^\ast$ denote the image of $V$ in $E^\ast$ and
note that $V^\ast$ is a $3${--}ball.
If $E^{S^1}=\R^2$ then $E^{S^1}\cap V^\ast$ is a $2${--}ball:
if $E^{S^1}$ is a point, then $E^{S^1}\cap V^\ast$ is a point
in the interior of the $3${--}ball.

In case $E^{S^1}=\R^2$, choose a smoothly embedded, 
compact, closed surface $Q\subset E^\ast$, 
so that $Q\,\cap\,\partial E^\ast$ is a $2${--}ball 
which contains both $V^\ast\,\cap\,\partial E^\ast$ and
$U^\ast\,\cap\,\partial E^\ast$ and so that the compact
component of $E^\ast-Q$ contains both $U^\ast$ and $V^\ast$.
In case $E^{S^1}$ is a point, the image of the singular set
has the following description.
There is one point, for the image of $E^{S^1}$,
together with a finite number
of proper rays, the images of the various $E^{K(p_i)}$'s.
Each ray crosses $\partial V^\ast=S^2$ transversely in one
point.
Choose $Q\subset E^\ast$ to be a smoothly embedded, 
compact, closed surface so that the compact
component of $E^\ast-Q$ contains both $U^\ast$ and $V^\ast$.
Further require that $Q$ intersects each of the rays transversely
in a single point.

Now write $E^\ast=N\cup W$ where $\partial N=\partial W=Q$
and $N$ is compact.
It follows from the Mayer{--}Vietoris theorem and intersection
theory that
$H_1(Q;\Z)=H_1(N;\Z)\oplus H_1(W;\Z)$ and that there is a symplectic
basis for $H_1(Q;\Z)$: $x_1$, \dots, $x_g$, $y_1$, \dots, $y_g$,
so that $x_i\cap x_j=y_i\cap y_j=0$ and $x_i\cap y_j=\delta_{ij}$.
Furthermore, the $x_i$ generate $H_1(N;\Z)$ and the $y_i$
generate $H_1(W;\Z)$ under the decomposition.

Represent the elements $y_i$ by disjoint embedded circles in $Q$:
in case $E^{S^1}$ is a point, arrange for these circles to miss the
points where the rays cross.
Let $\hat N$ be the result of doing surgery on these circles.
Each of these circles lies in the image of the principal orbit,
so the circle action over them is trivial.
Hence we can construct a compact, smooth $4${--}manifold $X^4$
with boundary which supports a smooth $S^1$ action and so that
$\pi^{-1}(N)$ is a smooth, equivariant submanifold of $X^4$.
The orbit space of the $S^1$ action on $X$ is just $\hat N$.
Let $J=\Bar{X-V}$ and let $J^\ast$ be the image of $J$ in $\hat N$.
Let $P$ denote the open dense subset of principal orbits in $J$.
If $E^{S^1}=\R^2$, $P=J-A$ where $A=J^{S^1}$ is an annulus.
If $E^{S^1}$ is a point, $P=X-\disjointunion_{i=1}^r A_i$, where 
each $A_i=J^{K(p_i)}$ is an annulus.

Consider the pair $(P,\partial V\cap P)$.
The $S^1$ action on this pair is free, and for the orbit space pair
$H^\ast(P^\ast,\partial V^\ast\cap P^\ast;\Z)=0$.
Hence $H^\ast(P,\partial V\cap P;\Z)=0$ by a spectral sequence argument.
From the description of $P$ in the last paragraph and 
the Mayer{--}Vietoris sequence $H^\ast(J,\partial V;\Z)=0$.
Since $V$ is a $4${--}ball, $X=V\cup J$ is
an integral homology $4${--}ball.
The double of $X$ is the required integral homology $4${--}sphere.\endprf

\newsec{A construction of $\R^4$'s with barrier $S^3$'s.}
In [\bT] we defined an invariant of smooth $\R^4$'s,
$\gamma$, which takes on integer values greater than or
equal to $0$ and $+\infty$.
We can only prove the existence of barriers in certain smoothings:
we call a smoothing {\sl definite\/} provided it is diffeomorphic at
$\infty$ to the end of some smoothing $M^4-pt$ where $M^4$ is a 
simply-connected, compact topological manifold 
with a definite intersection form
which can not be diagonalized over the integers.
In [\bT; eg~5.6] we constructed smoothings of $\R^4$, $E_n$, 
which are definite and which satisfy $\gamma(E_n)=n$, $0<n\leq\infty$:
indeed for each $n$, $0<n\leq\infty$, we construct a one parameter
family of them.

\newT{bst}
\proclaim\fbst
Let $E$ be a definite smoothing of $\R^4$ with $0<\gamma(E)<\infty$.
Then there exists a compact set $K\subset E$,
such that  any flat $S^3\subset E$
with $K$ on the inside of $S^3$ is a barrier.
Furthermore, $\gamma(inside)=\gamma(E)$
and the $inside$ is definite.\endproc

\noindent{\bf Remarks}\qua
We know of no example of a smoothing $E$ with 
$0<\gamma(E)\leq\infty$ which is not definite.
There are examples of $E$ with $\gamma(E)=\infty$
which are definite but do not have barrier $S^3$'s,
for example the universal $\R^4$ of [\bFT].
Any $E$ which embeds in the standard $\R^4$ has no barriers.

\prf
From [\bT; Theorem 5.1] we see that if $E_0\subset E_1$ are
smoothings of $\R^4$, $\gamma(E_0)\leq\gamma(E_1)$.
It follows from the definition that if $E_0$ is definite
so is $E_1$.

Call a neighborhood $U$ of $S^3$ a {\sl 0-neighborhood\/} if
$U$ is open and $U-S^3$ has two components.
We label the component which intersects the inside of $E-S^3$
the {\sl inside\/} and the other component the {\sl outside}.
Given any neighborhood $U$ of a flat $S^3$ in $E^4$,
we can find a smaller neighborhood homeomorphic to
$S^3\times(-\infty,\infty)$ and this is a 0-neighborhood.
Hence, to prove $S^3$ is a barrier, it suffices to prove
$e(S^3)\cap S^3\neq\emptyset$ for all smooth embeddings
$e\colon U\to E$ where $U$ is a 0-neighborhood.

Assume we have a flat $S^3\subset E$ and
let $U$ be a 0-neighborhood of $S^3$.
Let $e\colon U\to E_1$ be a smooth embedding of $U$
into any smoothing of $\R^4$.
Note $e(U)$ is a 0-neighborhood of $e(S^3)$.
\newL{iso}
\proclaim\fiso
Let $I$ denote the smoothing on the inside of $E-S^3$.
If $\gamma(I)>0$, $e$ takes the inside of $U$ to the inside of $e(U)$.
\endproc

\prf
If not, one can construct a smooth homotopy $4${--}sphere
with $I$ smoothly embedded.
But this contradicts $\gamma(I)>0$.\endprf

Again assume we have a flat $S^3\subset E$ and
let $U$ be a 0-neighborhood of $S^3$.
Let $e_1$, $e_2\colon U\to E_1$ be smooth embeddings of $U$
into the same smoothing of $\R^4$.
\newL{dso}
\proclaim\fdso
Let $I$ again denote the smoothing on the inside of $E-S^3$.
 Suppose $I$ is definite and $e_1(S^3)\cap e_2(S^3)=\emptyset$.
Then the inside of $e_1(S^3)$ and the inside of $e_2(S^3)$
are disjoint.\endproc

\prf
Note $\gamma(I)>0$ for any definite $\R^4$, so the ``inside''
is well{--}defined by \fiso.
If 3.3 were false, then we could construct a 
new smoothing of $M^4-pt$
with a periodic end in the sense of Taubes, [\bCT].
But this is precisely what the main theorem of [\bCT]
forbids.\endprf

We now return to the proof of the existence of barriers.
It follows from [\bT] that we can find an $S^3\subset E$
so that $\gamma(I)=n$ and $I$ is definite, where $I$ is the
inside of $E-S^3$.
We can further assume that $S^3\subset E$ has a smooth point.
Let $K$ denote the closure of $I$ in $E$.
Of course it is homeomorphic to a $4${--}ball.
Since $I$ is definite, it follows from [\bT; Theorem 5.3]
that we can find an integer $N>0$ such that
$\gamma(\natural^{2^{N}} I)>\gamma(E)$, where
$\natural^{2^{N}} I$ denotes the end-connected sum of $I$ with itself
$2^N$ times.

The proof proceeds by constructing successively larger
compact sets until the conclusion of the theorem holds.
We introduce some notation for the proof.
First we label the embedding
$S^3\subset E$ by $e^{(0)}_0\colon S^3\subset E$.
Then we label $I$ as $I^{(0)}$ and $K$ as $K^{(0)}$.
We will construct a sequence of embeddings,
$e^{(j)}_0\colon S^3\subset E$ starting with the $j=0$ we
have just exhibited.
Let $I^{(j)}$ be the inside of $E-e^{(j)}_0(S^3)$ and let
$K^{(j)}$ denote the closure of $I^{(j)}$.
As part of the construction, we will have
$K^{(j-1)}\subset I^{(j)}$.
We will continue the construction until the conclusion of
the theorem holds for $K^{(j)}$: this must happen for
some $j< N$ as we shall see.
Since $I^{(0)}$ is definite, 
$I^{(j)}$ is definite, and $0<\gamma(I^{(j)})=\gamma(E)<\infty$. 
A second part of the construction guarantees that
$\gamma(\natural^{2^{N-j}} I^{(j)})>\gamma(E)$.

\medskip\noindent$\bullet$\qua
{\sl Suppose we have constructed $e^{(j-1)}_0$.}\medskip

If the conclusion of the theorem holds with $K=K^{(j-1)}$ we are done.
If not, there exists an $e_1\colon S^3\subset E$ 
with $K^{(j-1)}$ on the inside and 
a $0$-neighborhood of $e_1(S^3)$, say $U$, which 
we may take to miss $K^{(j-1)}$,
so that there exists a smooth embedding 
$e\colon U\to E$
such that $e_1(S^3)\cap e(S^3)=\emptyset$.
By Lemmas 3.2 and 3.3, the insides
of $E-e_1(S^3)$ and $E-e(S^3)$ are disjoint.

Take a flat $S^3$ whose inside contains both $e_1(S^3)$ and 
$e(S^3)$
and whose boundary has a smooth point.
Denote the embedding by $e^{(j)}_0\colon S^3\subset E$ and
let $I^{(j)}$ denote the inside of $E-e^{(j)}_0(S^3)$.
By the definition of $\gamma$, $I^{(j)}$ smoothly embeds in $X^4$,
$X^4$ a closed, smooth, compact, Spin $4${--}manifold whose
rational intersection form is hyperbolic with $\gamma(E)$
hyperbolic summands.
Let $Y^4$ denote the result of removing the inside of $E-e(S^3)$
from $X$
and replacing it by the inside of $E-e_1(S^3)$:
$Y$ remains a closed, smooth, compact, Spin $4${--}manifold whose
rational intersection form is hyperbolic with $\gamma(E)$
hyperbolic summands and two disjoint copies of $I^{(j-1)}$ embed
smoothly in it, each copy having a smooth boundary point.
If $N-(j-1)=1$, then $\gamma(\natural^2 I^{(j-1)})>\gamma(I)$ 
and this contradicitions the embedding of two copies of
$I^{(j-1)}$ in $Y$.
Hence if $N-(j-1)=1$ the conclusion of the theorem must have held
and we are done.

If $N-(j-1)>1$, note the following.
By construction $I^{(j)}$ has 2 copies of $I^{(j-1)}$ embedded in it, 
each with a smooth boundary point.
Hence 
$\natural^{2^{N-(j-1)}} I^{(j-1)}\subset
\natural^{2^{N-j}} I^{(j)}$, so
$\gamma(\natural^{2^{N-j}} I^{(j)})>\gamma(E)$.

Now repeat the argument starting at $\bullet$ above with $j=j+1$.
Eventually $N-(j-1)=1$ and the process halts 
if it has not halted earlier.\endprf

\rk{Acknowledgement} The author is partially supported by the N.S.F.

\references
\Addresses\recd
\bye